\documentclass{amsart}

\usepackage{amssymb}
\usepackage[all]{xy}

\newtheorem{thm}{Theorem}%[section]

\newtheorem{oldthm}[thm]{Old Theorem}
\newtheorem{newthm}[thm]{New Theorem}

\newtheorem{prop}[thm]{Proposition}

\theoremstyle{definition}
\newtheorem{defn}[thm]{Definition}

\newtheorem{say}[thm]{}
\newtheorem{exmp}[thm]{Example}

\newtheorem{defn-thm}[thm]{Definition--Theorem}  %!!!!!!!!!!!!!!!!!!!!!!!!
\newtheorem{defn-lem}[thm]{Definition--Lemma}  %!!!!!!!!!!!!!!!!!!!!!!!!
\theoremstyle{remark}

\renewcommand{\o}[0]{{\mathcal O}} 
\newcommand{\z}[0]{{\mathbb Z}}

\renewcommand{\r}[0]{{\mathbb R}} 

\renewcommand{\a}[0]{{\mathbb A}}

\newcommand{\p}[0]{{\mathbb P}}

\newcommand{\q}[0]{{\mathbb Q}}
\newcommand{\map}[0]{\dasharrow}
\newcommand{\qtq}[1]{\quad\mbox{#1}\quad}

\newcommand{\supp}[0]{\operatorname{Supp}}    
\newcommand{\red}[0]{\operatorname{red}}    
\newcommand{\codim}[0]{\operatorname{codim}}

\newcommand{\chr}[0]{\operatorname{char}}

\newcommand{\rdown}[1]{\lfloor{#1}\rfloor}

\newcommand{\depth}[0]{\operatorname{depth}} 
\newcommand{\tsum}[0]{\textstyle{\sum}}

\def\loccoh#1.#2.#3.#4.{H^{#1}_{#2}(#3,#4)}

\DeclareMathAlphabet{\mathchanc}{OT1}{pzc}%
                                {m}{it}

\usepackage[all]{xy}\xyoption{dvips}

\newcommand{\vol}[0]{\operatorname{vol}}

\DeclareMathOperator{\bdp}{{\bf B}_{+}^{\rm div}}

\begin{document}
\bibliographystyle{amsalpha}

%\today

\title[Hilbert functions]{How much of the Hilbert function do we \\ really need to know?}
\author{J\'anos Koll\'ar}
\begin{abstract}The aim of this lecture is to describe several examples where
 the  leading coefficient of a Hilbert function
tells us everything we need.
\end{abstract}

\maketitle

The starting point is the following theorem, whose proof---though not
its precise statement---is in \cite[III.9.9]{hartsh}. 

\begin{oldthm} Let $f:X\to S$ be a projective morphism and $F$ a coherent sheaf on $X$. Then
\begin{enumerate}
\item $s\mapsto \chi\bigl(X_s, F_s(m)\bigr)$ is a lower semicontinuous function on $S$ for $m\gg 1$.
\item  If $S$ is connected and reduced, then  $F$ is flat over $S$  $\Leftrightarrow$ 
the above function  $s\mapsto \chi\bigl(X_s, F_s(m)\bigr)$ is   constant  on $S$ for every $m$.
\end{enumerate}
\end{oldthm}

Thus one can establish flatness
by computing the Hilbert function of the individual fibers $F_s$.
Note that  the fibers over points carry no information about the nilpotent directions in the base, so the restriction to reduced $S$ is necessary in (2).

In practice it is frequently quite hard to determine
the {\em whole} Hilbert function  $\chi\bigl(X, F(m)\bigr)$ for a coherent
sheaf $F$ on a proper scheme $X$, but it turns out that
 there are many interesting situations where 
it is enough to know the  leading coefficient of
 $\chi\bigl(X, F(m)\bigr)$ to guarantee flatness. The first such general result
I know of is  due to Hironaka \cite{MR0102519}; 
%% Ill J 1958
see also \cite[III.9.11]{hartsh}.
The projective case of the theorem can be formulated
as follows.

\begin{oldthm} 
Let $T$ be a connected, regular, 1-dimensional scheme and
$X\subset \p^N_T$ a closed subscheme, flat over $T$. Then
\begin{enumerate}
\item $t\mapsto  \deg (\red X_t)$ is a lower semicontinuous function on $T$.
\item If the reduced fibers
$\red X_t$ are  normal then  the following are equivalent.
\begin{enumerate}
\item $t\mapsto \deg (\red X_t)$ is    constant  on $T$,
\item $t\mapsto \chi\bigl(\red X_t, \o_{\red X_t}(m)\bigr)$ is   constant 
for every $m$ and
\item the fibers $X_t$ are reduced.
\end{enumerate}
\end{enumerate}
\end{oldthm}

The leading coefficient of $\chi\bigl(\red X_t, \o_{\red X_t}(m)\bigr)$
equals $\deg (\red X_t)/(\dim X_t)!$, thus we can informally
summarize the above theorem by saying that ``the leading coefficient
determines flatness.''

We are looking for theorems of this type. The first part should be a
general assertion that some invariants related to Hilbert functions are 
lower or upper  semicontinuous on the base. Then, under some geometric
assumptions, we aim to show that  constancy of the leading 
coefficient---usually given as the {\it volume} of a divisor as in
(\ref{hilb.funct.defn}.1)---implies
  constancy of the whole Hilbert function, hence flatness.

Each of the next 5 sections outlines such results.
A detailed treatment of the  claims in Section 1--2
will appear in \cite{k-modbook}. Sections 3--4 summarize some of the 
theorems of \cite{k-gl1, bha-dej, k-gl2} while
Section 5 is taken from \cite{fkl-vol}.   

\section{Simultaneous canonical models}

There are, unfortunately,  two distinct definitions of canonical models
in use.

\begin{defn}[Canonical models]\label{2.canmod.defns}
Let $(X, \Delta)$  be a proper log canonical pair such that
 $K_X+\Delta$ is big.  As in
 \cite[3.50]{km-book}, 
its {\it canonical model}  is the unique log canonical
pair   $(X^c, \Delta^c)$ such that $K_{X^c}+\Delta^c$ is ample
and 
$$
\tsum_{m\geq 0} H^0\bigl(X, \o_X(mK_X+\rdown{m\Delta})\bigr)\cong
\tsum_{m\geq 0} H^0\bigl(X^c, \o_{X^c}(mK_{X^c}+\rdown{m\Delta^c})\bigr).
$$
There is a natural birational map
$$
\phi: (X, \Delta) \map (X^c, \Delta^c).
\eqno{(\ref{2.canmod.defns}.1)}
$$
On the other hand, if $X$ is a proper variety with arbitrary singularities, 
then one can take a resolution  $X^r\to X$ and its
canonical model  $(X^r)^c$. Since this is independent of the choice of $X^r$, 
it is frequently called  the {\it canonical model}  of $X$.
I suggest to call it the {\it canonical model of resolutions} of $X$
and denote by  $X^{\rm cr}$.
 More generally, let $X$ be a proper, pure dimensional scheme over a field.
Start with  any resolution $X^r\to \red X$ and let  $X^{\rm cr}$
denote  the  disjoint union of the
canonical models of those components that are of general type.
With a slight abuse of terminolgy, there is a natural  map
$$
\phi: X\map X^{\rm cr},
\eqno{(\ref{2.canmod.defns}.2)}
$$
which is birational on the general type components and not
defined on the others.

If $X$ has log canonical singularities then both variants are defined. Note that
$X^{\rm c}\cong X^{\rm cr}$  if  $X$ has only canonical  singularities 
but not in general.
\end{defn}

\begin{defn}[Simultaneous canonical model]\label{sim.canmodel.defn}
 Let $f:X\to S$ be a proper  morphism
of pure relative dimension $n$. One can define 
a {\it simultaneous canonical model of resolutions} 
 $f^{\rm scr}:X^{\rm scr}\to S$. If we also have a divisor $\Delta$ on $X$ such that 
 the fibers  $(X_s,  \Delta_s)$ are log canonical,
then one can  also define a {\it simultaneous canonical model} 
 $f^{\rm sc}:\bigl(X^{\rm sc},\Delta^{\rm sc}\bigr) \to S$. These are given
by diagrams 
$$
\begin{array}{rcl} X &\stackrel{\phi}{\map} &X^{\rm scr}\\ f
&\searrow \quad\swarrow & f^{\rm scr}\\ 
 &S&
\end{array}
\qtq{respectively}
\begin{array}{rcl} (X,\Delta) &\stackrel{\phi}{\map} &
\bigl(X^{\rm sc},\Delta^{\rm sc}\bigr)\\ 
f&\searrow \quad\swarrow & f^{\rm sc}\\ 
 &S&
\end{array}
$$
where  $f^{\rm scr}$ and $f^{\rm sc}$ are flat, proper 
 and each 
$\phi_s:X_s
\map  X_s^{\rm scr}$ 
is  the  canonical model of the resolutions of $X_s$
(resp.\ each 
$\phi_s:(X_s,\Delta_s)
\map  \bigl(X_s^{\rm sc},\Delta_s^{\rm sc}\bigr) $ 
is  the  canonical model of $(X_s,  \Delta_s)$).

% Similarly,  let $f:(X,\Delta)\to S$ be a
%  morphism whose fibers  $(X_s,  \Delta_s)$ are log canonical.
%  A {\it simultaneous canonical model} 
% is a diagram 
% $$
% \begin{array}{rcl} (X,\Delta) &\stackrel{\phi}{\map} &
% \bigl(X^{\rm sc},\Delta^{\rm sc}\bigr)\\ 
% f&\searrow \quad\swarrow & f^{\rm sc}\\ 
%  &S&
% \end{array}
% \eqno{(\ref{sim.canmodel.defn}.2)}
% $$
% where  $f^{\rm sc}:X^{\rm sc}\to S$ is flat, proper 
%  and each 
% $\phi_s:(X_s,\Delta_s)
% \map  \bigl(X_s^{\rm sc},\Delta_s^{\rm sc}\bigr) $ 
% is  the  canonical model of $(X_s,  \Delta_s)$.
% \medskip

{\it Note:} % \ref{sim.canmodel.defn}.3. 
 We need the additional assumption that 
$f^{\rm sc}: \bigl(X^{\rm sc},\Delta^{\rm sc}\bigr)\to S$ be
 {\it locally stable,} equivalently, that 
$K_{X^{\rm sc}/S}+\Delta$ be $\q$-Cartier. See \cite{k-modsurv, k-modbook}  for
 discussions about this condition.
If the fibers $ X_s^{\rm sc}$ have canonical singularities then
$K_{X^{\rm sc}/S}$ is automatically $\q$-Cartier, thus we did not need to
assume this for simultaneous canonical models of resolutions.
\end{defn}

\begin{thm}[Numerical criterion for simultaneous canonical models I]
 \label{sim.canmodel.thm.0.1}
Let $S$ be a connected, seminormal  scheme of $\chr 0$
and $f:X\to S$ a morphism of pure relative dimension $n$. Then 
\begin{enumerate}
\item  $s\mapsto \vol(K_{X_s^r})$
is a lower semicontinuous function on $S$ and
\item $f: X \to S$ has a simultaneous canonical model of resolutions iff
this function is  constant (and positive).
\end{enumerate}
\end{thm}

Part (1)  was  first observed and proved
in \cite{MR837624, MR946250}.
% \item General theory predicts that a 
% simultaneous canonical model of resolutions excists iff the  Hilbert function  $\chi\bigl(X^{rc}_s, \o(mK_{X^{rc}_s})\bigr)$
% is locally constant for every $m$. 
% \item Usually it is easy to generalize similar proofs
% from smooth varieties to klt or lc pairs, but here
% adding any boundary can ruin the argument and the conclusion
%  as show by the
% Examples \ref{no.simcanmod.exmp.0}--\ref{mmp.does.not.restrict.exmp}.
% \end{itemize}

The following is  a similar result for normal 
lc pairs, but the   lower semicontinuity of
Theorem \ref{sim.canmodel.thm.0.1} changes to upper semicontinuity.

\begin{thm}[Numerical criterion for simultaneous canonical models II]
 \label{sim.canmodel.thm.0}
 Let
$S$ be a connected, seminormal scheme of $\chr 0$ and $f:(X,\Delta)\to S$ a
flat morphism whose fibers  $(X_s,  \Delta_s)$ are log canonical.
Then 
\begin{enumerate}%\setcounter{enumi}{3}
\item  $s\mapsto \vol(K_{ X_s}+\Delta_s)$
is an upper semicontinuous function on $S$ and
\item $f:(X, \Delta)\to S$ has a simultaneous canonical model iff
this function is  constant.
\end{enumerate} 
\end{thm}

{\it Notes.} Strictly speaking, part (2) needs the assumption that the fibers
$(X_s,  \Delta_s)$ have a canonical model. This is conjectured to
be true and it is known in many cases, for instance when $(X_s,  \Delta_s)$
is klt.

A stronger version of the theorem assumes only that each fiber is
normal in codimension 1 and has log canonical normalization.

A key ingredient of the proof of
Theorems \ref{sim.canmodel.thm.0.1}--\ref{sim.canmodel.thm.0}
is the following characterization of 
canonical models. 

\begin{prop} \label{canmod.vol.minimizer}
Let $X$ be a smooth proper variety of dimension $n$. Let $Y$ be a 
normal, proper  variety
birational to $X$ and $D$ an effective $\q$-divisor on $Y$
such that $K_Y+D$ is $\q$-Cartier, nef and big. Then
\begin{enumerate}
\item $\vol(K_X)\leq \vol( K_Y+D )=( K_Y+D )^n$ and
\item equality holds iff $D=0$ and $Y$ has canonical singularities.
\end{enumerate}
\end{prop}

\section{Simultaneous canonical  modifications}

For surfaces, the existence criterion of simultaneous canonical modifications
 is proved in \cite[Sec.2]{ksb}.
%The previous {\em upper} semicontinutity 
%condition changes to {\em lower} semicontinutity. 
In higher dimensions  we need  to work with a
sequence of intersection numbers and with their lexicographic ordering.

\begin{defn} \label{I(A,B).defn}
Let $X$ be a proper scheme of dimension $n$ and
$A,B$ $\r$-Cartier divisors on $X$. Their 
{\it sequence of intersection numbers} is
$$
I(A,B):=\bigl( (A^n),  \dots, 
(A^{n-i}\cdot B^i),
\dots, (B^n)\bigr)\in \r^{n+1}.
%\eqno{(\ref{lex.upp.sc.defn}.3)}
$$
For two divisors, the relevant Hilbert function is the 2-variable polynomial
$h(u,v):=\chi\bigl(X, \o_X(uA+vB)\bigr)$ and the above intersection numbers
are the coefficients of the leading homogeneous term, which has degree
$=\dim X$. 

The {\it lexicographic} ordering  is denoted by
$(a_0,\dots, a_n)\preceq (b_0,\dots, b_n)$.
(This holds if either $a_i=b_i$ for every $i$ or there is
an $r\leq n$ such that $a_i=b_i$ for $i<r$ but $a_r<b_r$.)
For polynomials we define an ordering
$$
f(t)\preceq g(t)\ \Leftrightarrow\  f(t)\leq g(t)\ \forall  t\gg 0.
$$
Note that
 $\tsum_i a_i t^{n-i}\preceq\tsum_i b_i t^{n-i} $ iff $
(a_0,\dots, a_n)\preceq (b_0,\dots, b_n)$.
Thus
%If we have proper schemes $X,X'$     of dimension $n$ and
%$\q$-Cartier divisors $A,B$  on $X$ and $A',B'$  on $X'$ then
$$
I(A,B)\preceq I(A',B')\ \Leftrightarrow\
(mA+B)^n \leq (mA'+B')^n\quad \forall  m\gg 0.
$$
\end{defn}

\begin{defn}[Simultaneous canonical modification]\label{sim.canmod.defn}

Let $Y$ be a scheme over a field $k$. 
(We allow $Y$ to be reducible and nonreduced.)
Its {\it canonical modification}  is a morphism $p:Y^{\rm can}\to Y$
such that  
$Y^{\rm can}\to \red Y$ is proper, birational, $Y^{\rm can}$ has  canonical  singularities and
$K_{Y^{\rm can}}$ is ample over $ Y$.

Let $\Delta$ be an effective divisor on $Y$.  A
{\it canonical modification}  is a morphism
 $p:\bigl(Y^{\rm can}, \Delta^{\rm can}\bigr)\to 
(Y,\Delta)$ where $p$ is  proper, birational,
$\Delta^{\rm can}=p^{-1}_*\Delta $, 
$\bigl(Y^{\rm can}, \Delta^{\rm can}\bigr) $ is canonical and
$K_{Y^{\rm can}}+\Delta^{\rm can}$ is ample over $ Y$. A canonical modification
is unique and it exists iff the following conditions hold:
\begin{enumerate}
\item[$(\diamond)$]  The reduced scheme $\red Y$ is
 smooth at the generic points of $\supp\Delta$
and all coefficients in $\Delta$ are in the interval $[0,1]$.
\end{enumerate}

  Let $f:X\to S$ be a morphism 
 of pure relative dimension $n$ and $\Delta$  an
effective  divisor on $Y$.
A {\it simultaneous canonical modification}
is a proper morphism $p:\bigl(X^{\rm scan}, \Delta^{\rm scan}\bigr)\to (X,\Delta)$ such that
 $f\circ p:\bigl(X^{\rm scan}, \Delta^{\rm scan}\bigr)\to S$ is locally stable (\ref{sim.canmodel.defn}.3) 
and
$p_s:\bigl(X^{\rm scan}_s, \Delta^{\rm scan}_s\bigr) 
\to (X_s,\Delta_s)$ 
is  the  canonical modification  for every $s\in S$.
\end{defn}

Let $S$ be a connected, seminormal scheme of $\chr 0$,
$f:X\to S$ a morphism of pure relative dimension $n$, $H$
an $f$-ample divisor class  and
$\Delta$ an effective divisor on $X$ such that
$(X_s,\Delta_s)$  satisfies the assumptions $(\diamond)$ %(\ref{sim.canmod.defn}.1)
for every $s\in S$. Thus  the canonical modifications 
$p_s:\bigl(X_s^{\rm can}, \Delta_s^{\rm can}\bigr)\to (X_s,\Delta_s)$ 
exist.

\begin{thm}[Numerical criterion for simultaneous canonical modification]
 \label{sim.canmod.thm.0}
With the above notation, 
\begin{enumerate}
\item  $s\mapsto I\bigl(p_s^*H_s, K_{X_s^{\rm can}}+\Delta_s^{\rm can}\bigr)$
is a lexicographically lower semicontinuous function on $S$ and
\item $f:(X,\Delta)\to S$ has a simultaneous canonical modification iff
this function is  constant.
\end{enumerate}
\end{thm}

There is also a similar condition for
simultaneous log canonical and semi-log-canonical modifications
but these only apply when $K_{X/S}+ \Delta $ is $\q$-Cartier.
The following example illustrates the problems that occur in general.

\begin{exmp} \label{no.simult.lcmod.exmp.1}
In $\p^2$ consider  a line  $L\subset \p^2$ and
a family of degree 8 curves $C_t$ such that
 $C_0$ has 4 nodes on $L$ plus an ordinary 6-fold point outside $L$ and
 $C_t$ is smooth and tangent to $L$ at 4 points for $t\neq 0$.

Let $\pi_t:S_t\to \p^2$ denote the double cover of $\p^2$ ramified along $C_t$.
Note that $K_{S_t}=\pi_t^*\o(1)$, thus $(K_{S_t}^2)=2$.
For each $t$, the preimage $\pi_t^{-1}(L)$ is a union of 2 curves
$D_t+D'_t$. Our example is the family of  pairs
$(S_t, D_t)$. We claim that
\begin{enumerate}
\item there is a log canonical modification
$ \bigl(S_t^{\rm lc}, D_t^{\rm lc}\bigr)  \to  (S_t, D_t)$  for every $t$ and
\item $\bigl(K_{S_t^{\rm lc}}+D_t^{\rm lc}\bigr)^2=1$  for every $t$ yet
\item there is no simultaneous log canonical  modification.
\end{enumerate}

If $t\neq 0$ then $S_t$ is smooth and $D_t$ is smooth.
Furthermore $D_t, D'_t$ meet transversally at 4 points, thus
$(D_t\cdot D'_t)=4$.  Using   
$\bigl((D_t+D'_t)^2\bigr)=2$, we obtain
 that $( D_t^2)=-3$. Thus
$(K_{S_t}+D_t)^2=1$.

If $t= 0$ then $S_0$ is singular at 5 points. 
$D_0, D'_0$ meet transversally at 4 singular points of type $A_1$, thus
$(D_0\cdot D'_0)=2$. This gives that $( D_0^2)=-1$. Thus
$(K_{S_0}+D_0)^2=3$.
The pair $(S_0, D_0)$ is lc away from the preimage of the 6-fold point.
Let $q:T_0\to S_0$ denote the minimal resolution of this point.
The exceptional curve $E$ is smooth, has genus 2 and $(E^2)=-2$.
Thus $K_{T_0}=q^*K_{S_0}-2E$ hence  $(T_0, E+D_0)$ is the
log canonical modification of $(S_0, D_0)$ and
$$
(K_{T_0}+E+D_0)^2=\bigl(q^*K_{S_0}-E+D_0\bigr)^2=
(K_{S_0}+D_0)^2+(E^2)=1.
$$
Thus  $\bigl(K_{S_t^{\rm lc}}+D_t^{\rm lc}\bigr)^2=1$ for every $t$.

Nonetheless, the log canonical  modifications do not form a flat family.
Indeed, such a family would be a family of surfaces with
ordinary nodes, so the relative canonical class would be a
Cartier divisor. However, $(K_{S_t}^2)=2$ for $t\neq 0$ but
$(K_{T_0}^2)=\bigl(q^*K_{S_0}-2E\bigr)^2=-6$.
\end{exmp}

\section{Families of Cartier divisors}

\begin{exmp} Consider the family of quadric surfaces
$$
X:=\bigl(x_1^2-x_2^2+x_3^2-t^2x_0^2=0\bigr)\subset \p^3_{\bf x}\times \a^1_t.
$$
The  fiber $X_0$ is a cone, the  other fibers  are smooth.
Consider the Weil divisors
$$
D:=\bigl(x_1-x_2=x_3-tx_0=0\bigr)\qtq{and}
E:=\bigl(x_1+x_2=x_3-tx_0=0\bigr).
$$
The fibers  $D_t, E_t$ form a pair of intersecting lines on $X_t$ for every $t$.
It is easy to compute that
% \begin{enumerate}
% \item $(aD_0+bE_0)^2=\tfrac12(a+b)^2$ and $aD_0+bE_0$ is Cartier iff
%  $a+b$ is even,
% \item $(aD_t+bE_t)^2=2ab$
% and $aD_t+bE_t$ is always Cartier for $t\neq 0$ and
% \item $aD+bE$ is Cartier iff $a=b$.
% \end{enumerate}
% From these we conclude that
\begin{enumerate}  %\setcounter{enumi}{3}
\item $(aD_0+bE_0)^2=\tfrac12(a+b)^2\geq 2ab=(aD_t+bE_t)^2$ and
\item equality holds iff $a=b$ iff $aD+bE$ is Cartier.
\end{enumerate}
\end{exmp}

We aim to prove that this example is quite typical, as far as intersection numbers are concerned. (It is, however, special in that the equations
define the restrictions $D_t$ unambiguously. In general, if $D$ is effective,
the sheaf theoretic restriction $\o_X(-D)|_{X_t}$ may have embedded points.
As long as the fibers are smooth in codimension 1, such embedded points
appear only in codimension $\geq 2$, so there is a well-defined Weil divisor
that can be thought of as the restriction $D_t$.) 

The following result was conjectured in 
\cite{k-gl1} and proved there for log canonical fibers. The extension to normal fibers is done in \cite{bha-dej}. 

\begin{thm}[Numerical criterion of  Cartier divisors, weak form]
\label{num.crit.glob.stab.0}
Let $C$ be a smooth, irreducible  curve
and $f:X\to C$ a proper, flat family of normal varieties of dimension $n$.
Let $D$ be a Weil divisor on $X$ such that its restriction $D_c$ is 
an ample Cartier divisor for every $c$.   Then 
\begin{enumerate}
\item  $c\mapsto  \bigl(D_c^n\bigr)$
 is an upper semicontinuous
 function on $C$ and
\item  $D$ is a Cartier divisor on $X$  iff
the above function  is   constant.
\end{enumerate}
\end{thm}

Ampleness is needed for $n\geq 3$, the main reason is that
$\bigl((-D)^n\bigr)=(-1)^n\bigl(D^n\bigr)$. Thus, on a 3--fold,
ample divisors behave anti-symmetrically while divisors pulled-back form a surface
behave symmetrically. 

The following general form is proved in 
\cite{k-gl2}, building on the earlier results of
\cite{k-gl1, bha-dej}.

 \begin{thm}[Numerical criterion of  Cartier divisors]
\label{num.C.crit.thm}
%Using the notation of (\ref{gen.numcrit.ass}) assume in addition that
Let $S$ be a connected, reduced scheme  over a field,
  $f:X\to S$  a flat, proper morphism of pure relative dimension $n$
with $S_2$ fibers
and $Z\subset X$ a closed subset 
such that  $\codim_{X_s}(Z\cap X_s)\geq 2$ for every $s\in S$.
Let  $L_U$ be an invertible sheaf on
$U:=X\setminus Z$ and assume that
the restriction $L_U|_{U_s}$ extends to an invertible sheaf $L_s$ on
$X_s$ for every $s\in S$.
 Then  
\begin{enumerate}
\item  $s\mapsto \bigl(H_s^{n-2}\cdot L_s^2)$
 is an upper semicontinuous
 function on $S$ and
\item  $L_U$ extends to an invertible sheaf $L$ on
$X$  iff the
 above function  is   constant.
\end{enumerate}
 Furthermore, if  $L_s$ is ample for every $s$ then  
\begin{enumerate}\setcounter{enumi}{2}
\item  $s\mapsto (L_s^n)$
 is an upper semicontinuous
 function on $S$ and
\item  $L_U$ extends to an invertible sheaf $L$ on
$X$  iff the
 above function  is   constant.
\end{enumerate}
\end{thm}

Note that taking  $\bigl(H_s^{n-2}\cdot\ $ in (1) is equivalent to
restricting to the intersection of $n-2$ very ample divisors.
In particular, the assumptions in (1) do not depend on singularities of the
fibers that appear in codimension $\geq 3$.
This is a  key point in the proof of Theorems \ref{num.crit.glob.stab.0}--\ref{num.C.crit.thm},
to be discussed next.

\section{Grothendieck--Lefschetz theorems for the local Picard group}

Let us recall the form given in \cite{sga2}.

\begin{oldthm}[Grothendieck--Lefschetz]\cite[XIII.2.1]{sga2}
Let  $(x\in X)$  be an excellent  local scheme, $x\in D\subset X$ a Cartier divisor.
Set  $U:=X\setminus\{x\}$, $U_D:=D\setminus\{x\}$ and let
   $L_U$ be a line bundle  on $U$ such that
    $L_U|_{U_D}\cong \o_{U_D}$.
\begin{enumerate}
\item[$(*)$]   Assume that  $\depth_x\o_D\geq 3$.
\end{enumerate}
Then    $L_U\cong \o_U$.
\end{oldthm}

For our pusposes, three aspects of this theorem are worth thinking about.

$\bullet$  It does not imply the usual Lefschetz theorem for hyperplane sections since a cone over a smooth projective variety
is usually only $S_2$ at its vertex.

$\bullet$  We would like to apply it to families of varieties
over a smooth curve $f:X\to C$ with $D$ being a fiber.  In this context
assuming that the fibers are $S_2$ is natural but $S_3$ is not. 
For instance, log canonical (and semi-log-canonical) varieties are
$S_2$ but frequently not $S_3$.

$\bullet$  The original form of the theorem assumes only that
$L$ is a rank 1 reflexive sheaf and in that setting the assumption
$(*)$ is optimal. However, in many potential applications
we know by induction that $L$ is  locally free on $U$. 
The following strengthening was conjectured in 
\cite{k-gl1} and proved there for log canonical fibers. The extension to normal fibers is done in \cite{bha-dej}, aside from some $p$-torsion
questions in characteristic $p$. The general form below
is established in \cite{k-gl2}.
Conjecturally, the result should hold for any excellent local scheme, but
the current proofs do not work in mixed characteristic.

\begin{newthm}
Let  $(x\in X)$  be a  local scheme that is essentially of finite type  over a field and
 $x\in D\subset X$ a Cartier divisor.
Set  $U:=X\setminus\{x\}$, $U_D:=D\setminus\{x\}$ and let
   $L_U$ be a line bundle  on $U$ such that
    $L_U|_{U_D}\cong \o_{U_D}$.
\begin{enumerate}
\item[$(**)$]   Assume that  $\depth_x\o_D\geq 2$ and $\dim_x D\geq 3$.
\end{enumerate}
Then    $L_U\cong \o_U$.
\end{newthm}

\begin{say}[Proof of the old form] 
Let $t$ be a defining equation of $D$ and write $L_D:=L_U|_{U_D}$. The sequence
 $0\to L_U\stackrel{t}{\to} L_U\stackrel{r}{\to} L_D\cong \o_{U_D}\to 0$ gives
$$
\begin{array}{cccccc}
H^0\bigl(U,L_U\bigr) &\stackrel{t}{\to} &
H^0\bigl(U,L_U\bigr)
&\stackrel{r}{\to} &H^0\bigl(U_D, L_D\cong \o_{U_D}\bigr) & \to\\
H^1\bigl(U, L_U\bigr) &\stackrel{t}{\to} &
H^1\bigl(U, L_U\bigr)
&\to &H^1\bigl(U_D,L_D\cong  \o_{U_D}\bigr).
\end{array}
$$
The assumption  $\depth_x\o_D\geq 3$ implies that
$H^1\bigl(U_D, \o_{U_D}\bigr)=0$ (see \cite[Sec.3]{gro-loc-coh-MR0224620}) and so the map
 $t:H^1\bigl(U, L_U\bigr) {\to} H^1\bigl(U, L_U\bigr)$ is surjective.
Next, $\dim U\geq 4$ implies  that $H^1\bigl(U, L_U\bigr)$ has finite length
(see \cite[VIII.2.3]{sga2}),
which implies that the map   $t:H^1\bigl(U, L_U\bigr) {\to} H^1\bigl(U, L_U\bigr)$ is an isomorphism.

Therefore $r: H^0\bigl(U,L_U\bigr)
{\to} H^0\bigl(U_D,L_D\bigr)$
is surjective and the 
constant 1 section of $L_D\cong \o_{U_D}$ lifts back to 
a nowhere-zero section of $L_U$. \qed
\end{say}

The vanishing $H^1\bigl(U_D, \o_{U_D}\bigr)=0$ is pretty much
equivalent to  $\depth_x\o_D\geq 3$, so the argument does not work if
 $\depth_x\o_D=2$.
  
Bhatt and de~Jong observed that one can go around this problem
in positive characteristic as follows. Assume that $X$ is normal
and let $X^+\to X$ denote the normalization of $X$ in an
algebraic closure of its field of functions $k(X)$.
Then $X^+$ is non-noetherian but it is CM by \cite{hoc-hun}. 
We can  lift everything back to $X^+$, apply the above proof
and then descend to $X$ at the end. There are several foundational
issues to deal with while working on $X^+$ (see \cite{bha-dej}) and the descent 
proves only that
$L_U^m\cong \o_U$ for some $m>0$. 

It is technically simpler to view $\o_{X^+}$ as a
quasi-coherent sheaf on $X$ and work with it; see \cite{k-modbook}.

Lifting back to characteristic 0
is easier. The extension to the non-normal case relies on the
structure theory of the local Picard group developed in
\cite{k-gl2}.

\section{Variation of $\r$-divisors}

This topic has the same spirit as the previous ones and
it is also used in the proofs of the theorems in Section 1--2.

\begin{defn}\label{hilb.funct.defn}
Let $X$ be a proper, normal algebraic variety of dimension $n$ over a 
field $K$ and $D$ an $\mathbb R$-divisor on $X$.
The {\it Hilbert function} of $D$ is the function
$$
{\mathcal H}(X,D): m\mapsto 
h^0( mD):=\dim_{K} H^0(X,\o_X(\rdown{ mD}));
%\eqno{(\ref{hilb.funct.defn}.1)}
$$
defined for all $m\in \r$.
If $D$ is an ample Cartier divisor then ${\mathcal H}(X,D) $
agrees with the usual Hilbert polynomial whenever $m\gg 1$
is an integer, but in general
 ${\mathcal H}(X,D)$ is not  a polynomial, not  even if $D$ is a $\z$-divisor and 
  $m\in \z$.
The simplest numerical invariant associated to the  Hilbert function
is the {\it volume} of $D$, defined as 
 $$\vol(D):=\limsup_{m\to\infty}\frac{ h^0(mD)}{m^n/n!}.
\eqno{(\ref{hilb.funct.defn}.1)}
$$
The volume is preserved by $\r$-linear
equivalence but the Hilbert function is not; see Example \ref{Hf.lineq.exmp}.
If $E$ is an effective $\r$-divisor, then clearly
$$
h^0( mD-mE)\leq h^0( mD)\leq h^0( mD+mE)
%\eqno{(\ref{hilb.funct.defn}.3)}
$$
holds for every $m>0$, hence
$\vol(D-E)\leq  \vol(D)\leq \vol(D+E)$.
\end{defn}

We claim that, although the volume does not determine the Hilbert function,
the only way to change the Hilbert function by subtracting or adding an effective divisor is to change the volume.

\begin{thm}\cite{fkl-vol}
 Let  $X$ be a proper, normal algebraic variety  over a perfect field, $D$ a big $\r$-divisor on $X$ and $E$ an effective $\r$-divisor on $X$. Then 

{\rm (Subtraction version.)} The following are equivalent.
\begin{enumerate}
\item[($1^-$)] $\vol(D-E)=\vol(D)$.
\item[($2^-$)] $h^0(mD-mE)=h^0(mD)$ for all $m>0$.
\item[($3^-$)] $E\leq N_{\sigma}(D)$, the negative part  of  the  Zariski--Nakayama-decomposition.
\end{enumerate}

{\rm (Addition version.)} The following are equivalent.
\begin{enumerate}
\item[($1^+$)] $\vol(D+E)=\vol(D)$.
\item[($2^+$)] $h^0(mD+mE)=h^0(mD)$ for all $m>0$.
\item[($3^+$)] $\supp(E)\subseteq\bdp(D)$, the divisorial part of the 
 augmented base locus of $D$.
\end{enumerate}
\end{thm}

\begin{exmp}\label{Hf.lineq.exmp}
Let $S\to \mathbb P^1$ be a minimal ruled surface with a negative section
$E\subset S$ and a positive section  $C\subset S$ that is disjoint from
$E$. Let $F_1,\dots, F_4$ be distinct fibers.
Then $C\sim_{\mathbb R} C+ (F_1-F_2)+\sqrt{2}(F_3-F_4).$

Note that
$\rdown{mC+ m(F_1-F_2)+m\sqrt{2}(F_3-F_4)} $
has negative intersection with $E$ for all real $m>0$. This implies that,
for every $m>0$ we have
$$
h^0\bigl(S, mC+ m(F_1-F_2)+m\sqrt{2}(F_3-F_4)\bigr)<
h^0\bigl(S, mC\bigr).
$$

\end{exmp}

%\bibliography{refs-main/refs}
\def\cprime{$'$} \def\cprime{$'$} \def\cprime{$'$} \def\cprime{$'$}
  \def\cprime{$'$} \def\cprime{$'$} \def\dbar{\leavevmode\hbox to
  0pt{\hskip.2ex \accent"16\hss}d} \def\cprime{$'$} \def\cprime{$'$}
  \def\polhk#1{\setbox0=\hbox{#1}{\ooalign{\hidewidth
  \lower1.5ex\hbox{`}\hidewidth\crcr\unhbox0}}} \def\cprime{$'$}
  \def\cprime{$'$} \def\cprime{$'$} \def\cprime{$'$}
  \def\polhk#1{\setbox0=\hbox{#1}{\ooalign{\hidewidth
  \lower1.5ex\hbox{`}\hidewidth\crcr\unhbox0}}} \def\cdprime{$''$}
  \def\cprime{$'$} \def\cprime{$'$} \def\cprime{$'$} \def\cprime{$'$}
\providecommand{\bysame}{\leavevmode\hbox to3em{\hrulefill}\thinspace}
\providecommand{\MR}{\relax\ifhmode\unskip\space\fi MR }
\providecommand{\MRhref}[2]{%
  \href{http://www.ams.org/mathscinet-getitem?mr=#1}{#2}
}
\providecommand{\href}[2]{#2}

\medskip

\noindent Princeton University, Princeton NJ 08544-1000

{\begin{verbatim} kollar@math.princeton.edu\end{verbatim}}

\end{document}